\documentclass[reqno,12pt]{amsart}
\usepackage{graphicx}
\usepackage[top=1in,bottom=1in,left=1in,right=1in]{geometry}
\newtheorem{thm}{Theorem}[section]
\newtheorem{cor}[thm]{Corollary}
\newtheorem{lem}[thm]{Lemma}
\newtheorem{prop}[thm]{Proposition}
\theoremstyle{definition}
\newtheorem{defn}{Definition}[section]
\theoremstyle{remark}
\newtheorem{rem}{Remark}[section]
\numberwithin{equation}{section}

\begin{document}
\title[Optimal conditions for elliptic systems]
      {Optimal conditions for $L^\infty$-regularity and a priori estimates for elliptic systems, I: two components}%

\author[Li Yuxiang]{Li Yuxiang }%
\address{Department of Mathematics, Southeast University, Nanjing 210096, P. R. China\\
         and\\
         Laboratoire Analyse, G\'{e}om\'{e}trie et Applications, Institut Galil\'{e}e,
         Universit\'{e} Paris-Nord 93430 Villetaneuse, France}
\email{lieyx@seu.edu.cn}
\thanks{Supported in part by National Natural Science Foundation of China 10601012
        and Southeast University Award Program for Outstanding Young Teachers 2005.}



\begin{abstract}
In this paper we present a new bootstrap procedure for elliptic
systems with two unknown functions. Combining with the
$L^p$-$L^q$-estimates, it yields the optimal $L^\infty$-regularity
conditions for the three well-known types of weak solutions:
$H_0^1$-solutions, $L^1$-solutions and $L^1_\delta$-solutions.
Thanks to the linear theory in $L^p_\delta(\Omega)$, it also yields
the optimal conditions for a priori estimates for
$L^1_\delta$-solutions. Based on the a priori estimates, we improve
known existence theorems for some classes of elliptic systems.
\end{abstract}
\maketitle

\section{Introduction}

The aim of this paper is to present a new alternate-bootstrap
procedure to obtain $L^\infty$-regularity and a priori estimates for
solutions of semilinear elliptic systems. This method enables us to
obtain the optimal $L^\infty$-regularity conditions for the three
well-known  types of weak solutions: $H_0^1$-solutions,
$L^1$-solutions and $L^1_\delta$-solutions of elliptic systems (for
their definitions, see Section 2). Combining with the linear theory
in $L^p_\delta$-spaces, our method also enables us to obtain a
priori estimates for $L^1_\delta$-solutions, therefore to improve
existence theorems for various classes of elliptic systems.

Let us consider the Dirichlet system of the form
\begin{eqnarray}
   &-\Delta u=f(x,u,v),&\ \ {\rm in}\ \Omega, \nonumber\\
   &-\Delta v=g(x,u,v),&\ \ {\rm in}\ \Omega, \label{sys:main}\\
   &u=v=0,&\ \ {\rm on}\ \partial\Omega,\nonumber
\end{eqnarray}
where $\Omega\subset \mathbb{R}^n$ is a smoothly bounded domain and
$f,g:\Omega\times \mathbb{R}^2\rightarrow \mathbb{R}$ are
Carath\'{e}odory functions. A typical case is
\begin{eqnarray}
   &-\Delta u=u^rv^p,&\ \ {\rm in}\ \Omega, \nonumber\\
   &-\Delta v=u^qv^s,&\ \ {\rm in}\ \Omega, \label{sys:secondary}\\
   &u=v=0,&\ \ {\rm on}\ \partial\Omega,\nonumber
\end{eqnarray}
where $r,s\geq 0,\ p,q>0$.

As a motivation, let us mention that in an important recent article
\cite{QS}, Quittner \& Souplet developed an alternate-bootstrap
method in the scale of weighted Lebesgue spaces
$L^p_\delta(\Omega)$. Their bootstrap procedure works well for
system (\ref{sys:main}) with
\begin{eqnarray}
   &&-h_1(x)\leq f\leq C_1(|v|^p+|u|^\gamma)+h_2(x),\nonumber\\
        [-1.5ex]&&\hskip 70mm u, v\in \mathbb{R},\ x\in\Omega,\label{ass:QS}\\[-1.5ex]
   &&-h_1(x)\leq g\leq C_1(|u|^q+|v|^\sigma)+h_2(x),\nonumber
\end{eqnarray}
where $p,q>0,\ pq>1,\ \gamma,\sigma\geq 1,\ C_1>0$, $h_1\in
L^1_\delta(\Omega)$, $h_2\in L^\theta$ with $\theta>n/2$.   They
obtained the optimal conditions for $L^\infty$-regularity and a
priori estimates for $L^1_\delta$-solutions, see \cite[Theorem
2.1]{QS}. The optimality was shown by Souplet \cite[Theorem 3.3]{S}.
Using this method, they obtained new existence theorems for various
classes of elliptic systems.

Our bootstrap procedure works for system (\ref{sys:main}) with $f,g$
satisfying more general assumptions
\begin{eqnarray}
   &&|f|\leq C_1(|u|^r|v|^p+|u|^\gamma)+h(x),\nonumber\\
        [-1.5ex]&&\hskip 60mm u, v\in \mathbb{R},\ x\in\Omega,\label{ass:f and g}\\[-1.5ex]
   &&|g|\leq C_1(|u|^q|v|^s+|v|^\sigma)+h(x),\nonumber
\end{eqnarray}
where $r,s,\gamma,\sigma\geq 0,\ p,q>0,\ C_1>0$ and the regularity
of $h$ will be specified later. The bootstrap procedure is only
based on the $L^m$-$L^k$-estimates in the linear theories of weak
solutions. So we are able to obtain the optimal
$L^\infty$-regularity conditions for the three well-known types of
weak solutions: $H_0^1$-solutions, $L^1$-solutions and
$L^1_\delta$-solutions of elliptic systems. Under some additional
appropriate conditions on $f,g$, this method also enables us to
obtain a priori estimates for $L^1_\delta$-solutions.

\subsection{Optimal conditions for $L^\infty$-regularity}
First we consider the case where $pq>(1-r)(1-s)$. Set
\begin{equation}\label{def:alpha,beta}
    \alpha=\displaystyle\frac{p+1-s}{pq-(1-r)(1-s)},\ \
    \beta=\displaystyle\frac{q+1-r}{pq-(1-r)(1-s)}.
\end{equation}
Note that $(\alpha,\beta)$ is the solution of
\begin{eqnarray*}
  \left[
  \begin{array}{cc}
    r-1 & p \\
    q & s-1
  \end{array}
  \right]
  \left[
  \begin{array}{c}
    \alpha \\
    \beta
  \end{array}
  \right]
  =
  \left[
  \begin{array}{cc}
    1 \\
    1
  \end{array}
  \right].
\end{eqnarray*}
Throughout this paper, we assume that $\alpha,\beta>0$, which is
obvious if $r,s\leq 1$.  The numbers $\alpha,\beta$ are related to
its scaling properties of system (\ref{sys:secondary}) (see for
instance \cite{CFMT}). For the parabolic counterpart of
(\ref{sys:secondary}), these numbers appear for instance in
\cite{DE, W, Zh} in the study of blow-up.

For the $L^\infty$-regularity, we obtain the following theorems.

\begin{thm}\label{thm:H1 solution}{\rm (Optimal $L^\infty$-regularity for
$H_0^1$-solutions)}\\
Assume that $f,g$ satisfy $(\ref{ass:f and g})$ with
$pq>(1-r)(1-s)$.
\begin{enumerate}
  \item[(i)] If
            \begin{eqnarray}\label{ass:optimal condition for H1}
                  &\displaystyle\max\{\alpha,\beta\}>\frac{n-2}{4},\ \ \
                  r,s,\gamma, \sigma<\frac{n+2}{n-2},\ \ \nonumber\\
                  [-1.5ex]\\[-1.5ex]
                  &\displaystyle\min\{p+r,q+s\}<\frac{n+2}{n-2},\ \
                  \ h\in L^\theta(\Omega),\
                  \theta>\frac{n}{2},\nonumber
           \end{eqnarray}
           then any $H_0^1$-solution of system $(\ref{sys:main})$ belongs to
           $L^\infty(\Omega)$;
  \item[(ii)]If $n\geq 3$ and
           \begin{eqnarray}\label{ass:optimal condition inverse for H1}
                 \displaystyle\max\{\alpha,\beta\}<\frac{n-2}{4},
           \end{eqnarray}
           system $(\ref{sys:main})$ in $B_1$, the unit ball in $\mathbb{R}^n$,
           with $f=(u+c_1)^r(v+c_2)^p$ and $g=(u+c_1)^q(v+c_2)^s$ for some $c_1,c_2>0$ admits a
           positive $H_0^1$-solution $(u,v)$ such that $u\notin
           L^\infty(B_1)$ and $v\notin L^\infty(B_1)$.
\end{enumerate}
\end{thm}

\begin{thm}\label{thm:L1 solution}{\rm (Optimal $L^\infty$-regularity for
$L^1$-solutions)}\\
Assume that $f,g$ satisfy $(\ref{ass:f and g})$ with
$pq>(1-r)(1-s)$.
\begin{enumerate}
  \item[(i)] If
           \begin{eqnarray}\label{ass:optimal condition for L1}
                  &\displaystyle\max\{\alpha,\beta\}>\frac{n-2}{2},\ \ \
                  r,s,\gamma, \sigma<\frac{n}{n-2},\ \ \nonumber\\
                  [-1.5ex]\\[-1.5ex]
                  &\displaystyle\min\{p+r,q+s\}<\frac{n}{n-2},\ \
                  \ h\in L^\theta(\Omega),\
                  \theta>\frac{n}{2},\nonumber
           \end{eqnarray}
             then any $L^1$-solution of system $(\ref{sys:main})$ belongs to
             $L^\infty(\Omega)$;
  \item[(ii)]If $n\geq 3$ and
             \begin{eqnarray}\label{ass:optimal condition inverse for L1}
                 \displaystyle\max\{\alpha,\beta\}<\frac{n-2}{2},
             \end{eqnarray}
             system $(\ref{sys:main})$ in $B_1$, the unit ball in $\mathbb{R}^n$,
             with $f=(u+c_1)^r(v+c_2)^p$ and $g=(u+c_1)^q(v+c_2)^s$ for some $c_1,c_2>0$ admits a
             positive $L^1$-solution $(u,v)$ such that $u\notin L^\infty(B_1)$ and $v\notin L^\infty(B_1)$.
\end{enumerate}
\end{thm}

\begin{thm}\label{thm:very weak solution}{\rm (Optimal $L^\infty$-regularity for
$L^1_\delta$-solutions)}\\
Assume that $f,g$ satisfy $(\ref{ass:f and g})$ with
$pq>(1-r)(1-s)$.
\begin{enumerate}
  \item[(i)] If
            \begin{eqnarray}\label{ass:optimal condition for very weak}
                  &\displaystyle\max\{\alpha,\beta\}>\frac{n-1}{2},\ \ \
                  r,s,\gamma, \sigma<\frac{n+1}{n-1},\ \ \nonumber\\
                  [-1.5ex]\\[-1.5ex]
                  &\displaystyle\min\{p+r,q+s\}<\frac{n+1}{n-1},\ \
                  \ h\in L^\theta_\delta(\Omega),\
                  \theta>\frac{n+1}{2},\nonumber
           \end{eqnarray}
            then any $L^1_\delta$-solution of system $(\ref{sys:main})$ belongs
            to $L^\infty(\Omega)$;
  \item[(ii)] If $n\geq 2$ and
            \begin{eqnarray}\label{ass:optimal condition inverse for very weak}
                 \displaystyle\max\{\alpha,\beta\}<\frac{n-1}{2},
            \end{eqnarray}
            there exist functions $a,b\in L^\infty(\Omega)$, $a,b\geq 0$ such
            that system $(\ref{sys:main})$ with $f=a(x)u^rv^p$ and $g=b(x)u^qv^s$ admits
            a positive $L^1_\delta$-solution $(u,v)$ such that $u\notin
            L^\infty(\Omega)$ and $v\notin L^\infty(\Omega)$.
\end{enumerate}
\end{thm}

Our theorems are closely related to the three critical exponents:
\begin{eqnarray*}
  &&p_S:=\left\{
  \begin{array}{ll}
    \infty & \mathrm{if}\ n\leq 2, \\
    (n+2)/(n-2) & \mathrm{if}\ n\geq 3,
  \end{array}
  \right.\\
  &&p_{sg}:=\left\{
  \begin{array}{ll}
    \infty & \mathrm{if}\ n\leq 2, \\
    n/(n-2) & \mathrm{if}\ n\geq 3,
  \end{array}
  \right.\\
  &&p_{BT}:=\left\{
  \begin{array}{ll}
    \infty & \mathrm{if}\ n\leq 1, \\
    (n+1)/(n-1) & \mathrm{if}\ n\geq 2.
  \end{array}
  \right.
\end{eqnarray*}
$p_S$ is the Sobolev exponent. $p_{sg}$ and $p_{BT}$ appear in study
of $L^1$-solutions and $L^1_\delta$-solutions of scalar elliptic
equations respectively. Note that
\begin{eqnarray*}
  \frac{n-2}{4}=\frac{1}{p_S-1},\ \frac{n-2}{2}=\frac{1}{p_{sg}-1},\
  \frac{n-1}{2}=\frac{1}{p_{BT}-1}.
\end{eqnarray*}
So if we write each critical exponent as $p_c$, the optimal
conditions for $L^\infty$-regularity of the above three types of
weak solutions have a consistent form
$\max\{\alpha,\beta\}>1/(p_c-1)$ and
$r,s,\gamma,\sigma,\min\{p+r,q+s\}<p_c$.

\begin{rem}
If $r,s\leq 1$, $\min\{p+r,q+s\}<p_c$ in Theorem \ref{thm:H1
solution}-\ref{thm:very weak solution} (i) is superfluous, see
Remark \ref{rem:p+r<pc}.
\end{rem}

For $pq\leq (1-r)(1-s)$, we have the following theorem.

\begin{thm}\label{thm:pq leq (1-r)(1-s)}
Assume that $f,g$ satisfy $(\ref{ass:f and g})$ with
$pq\leq(1-r)(1-s)$. Then Theorem \ref{thm:H1 solution}-\ref{thm:very
weak solution} $(i)$ also hold if $\max\{\alpha,\beta\}>1/(p_c-1)$
is replaced by $pq-(1-r)(1-s)<(p_c-1)\max\{p+1-s,q+1-r\}$.
\end{thm}

In order to justify the above theorems, let us recall the optimal
$L^\infty$-regularity for the scalar equation
\begin{eqnarray}
   &-\Delta u=f(x,u),&\ \ {\rm in}\ \Omega, \nonumber\\
     [-1.5ex]\label{eq:main}\\[-1.5ex]
   &u=0,&\ \ {\rm on}\ \partial\Omega,\nonumber
\end{eqnarray}
where $|f|\leq C(1+|u|^p)$ with $p\geq 1$. It is well-known that the
Sobolev exponent $p_S$ plays an important role in the optimal
$L^\infty$-regularity and a priori estimates of the
$H_0^1$-solutions, see \cite{FLN, GS, JL, ZZ} and the references
therein. Any $H_0^1$-solution of (\ref{eq:main}) belongs to
$L^\infty(\Omega)$ if and only if $p\leq p_S$, see for instance
\cite{BK, St}. For the $L^1$-solutions, the critical exponent is
$p_{sg}$. Any $L^1$-solution of (\ref{eq:main}) belongs to
$L^\infty(\Omega)$ if and only if $p<p_{sg}$, see for instance
\cite{A, NS, P}.

The critical exponent $p_{BT}$ first appeared in the work of
Br\'{e}zis \& Turner in \cite{BT}. They obtained a priori estimates
for all positive $H_0^1$-solutions of (\ref{eq:main}) for $p<p_{BT}$
using the method of Hardy-Sobolev inequalities. However the meaning
of $p_{BT}$ was clarified only recently. It was shown by Souplet
\cite[Theorem 3.1]{S} that $p_{BT}$ is the critical exponent for the
$L^\infty$-regularity of $L^1_\delta$-solutions of (\ref{eq:main})
by constructing an unbounded solution with $f=a(x)u^p$ for some
$a\in L^\infty(\Omega)$, $a\geq 0$ if $p>p_{BT}$. The critical case
$p=p_{BT}$ was recently shown to belong to the singular case for
$f=u^p$, see \cite{DMP}, also \cite{MR} for related results.
Moreover, the results of \cite{S} was extended to the case $f=u^p$
when $p>p_{BT}$ is close to $p_{BT}$.

If we set $\alpha=1/(p-1)$, i.e., the solution of $(p-1)\alpha=1$,
the optimal conditions for $L^\infty$-regularity of the above three
types of weak solutions also have a consistent form
$\alpha>1/(p_c-1)$. For more detailed discussions, we refer to the
book \cite[Chapter I]{QS$_2$}.

\vskip 3mm

Using the bootstrap procedure they developed based on linear theory
in $L^p_\delta(\Omega)$, Quittner \& Souplet \cite[Theorem 2.1]{QS}
obtained similar $L^\infty$-regularity condition as Theorem
\ref{thm:very weak solution} (i) assuming that $f,g$ satisfy
(\ref{ass:QS}). In \cite[Theorem 3.3]{S}, Souplet proved a similar
result as in Theorem \ref{thm:very weak solution} (ii) in the case
$f=a(x)v^p$ and $g=b(x)u^q$ for some functions $a,b\in
L^\infty(\Omega)$, $a,b\geq 0$.

\begin{rem}\label{rem:hyperbola}
Using the method of moving planes and Pohozaev-type identities, in
the case $f=v^p$ and $g=u^q$, $p,q>1$, it is proved if $\Omega$ is
convex and bounded, and $\alpha+\beta>(n-2)/2$, then there exists a
positive classical solution of (\ref{sys:main}); If $n\geq 3$,
$\Omega$ is starshaped and bounded, and $\alpha+\beta\leq (n-2)/2$,
then (\ref{sys:main}) has no positive solution, see \cite{CFM,
M$_2$}. Note that the optimal $L^\infty$-regularity condition in
Theorem \ref{thm:H1 solution} is weaker than the existence
condition, i.e., the so-called Sobolev hyperbola.
\end{rem}

\begin{rem}
We shall use a bootstrap procedure to prove the above theorems.
Based on another bootstrap procedure, using the method of
Rellich-Pohozaev identities and moving planes, \cite[Lemma 2.2]{CFM}
obtained a priori estimates for $H^1_0$-solutions of
(\ref{sys:main}) with $f,g$ satisfying some conditions similar to
(\ref{ass:QS}).
\end{rem}

\subsection{Optimal conditions for a priori estimates and existence theorems}

Combining with the linear theory in $L^p_\delta$-spaces, developed
in \cite{FSW}, see also \cite{BV}, our bootstrap procedure enables
us to obtain a priori estimates for system (\ref{sys:main}) with
$f,g$ satisfying (\ref{ass:f and g}) and
\begin{eqnarray}\label{ass:a priori estimates}
   f+g\geq-C_2(u+v)-h_1(x),\ \ \ u, v\in \mathbb{R},\ x\in\Omega,
\end{eqnarray}
where $C_2>0$, $h_1\in L^1_\delta(\Omega)$. By an a priori estimate,
we mean an estimate of the form
\begin{equation}\label{est:a priori estimates}
    \|u\|_{\infty}\leq C,\ \ \|v\|_{\infty}\leq C
\end{equation}
for all possible nonnegative solutions of (\ref{sys:main}) (in a
given set of functions), with some constant $C$ independent of
$(u,v)$. Our main result of the a priori estimates is the following
theorem.

\begin{thm}\label{thm:a priori estimates}
Let $f,g$ satisfy $(\ref{ass:f and g})$ and $(\ref{ass:a priori
estimates})$ with $pq>(1-r)(1-s)$ and $(\ref{ass:optimal condition
for very weak})$. Then there exists $C>0$ such that for any
nonnegative solution $(u,v)$ of $(\ref{sys:main})$ satisfying
\begin{equation}\label{ass:u v}
    \|u\|_{L^1_\delta}+\|v\|_{L^1_\delta}\leq M,
\end{equation}
it follows that $u,v\in L^\infty(\Omega)$ and
\begin{equation*}
    \|u\|_{L^\infty}+\|v\|_{L^\infty}\leq C.
\end{equation*}
The constant $C$ depends only on
$M,\Omega,p,q,r,s,\gamma,\sigma,C_1,C_2$.
\end{thm}

(\ref{ass:optimal condition for very weak}) is optimal for the a
priori estimates for the $L^1_\delta$-solutions of the system
(\ref{sys:main}) under the assumptions (\ref{ass:f and g}) and
(\ref{ass:a priori estimates}), see Theorem \ref{thm:very weak
solution} (ii).

There are several methods for the derivation of a priori estimates:
The method of Rellich-Pohozaev identities and moving planes, see
\cite{CFM, FLN}; The scaling or blow-up methods, which proceeds by
contradiction with some known Liouville-type theorems, see \cite{BM,
CFMT, FY, GS, Lou, So, Zou} and references therein, for the related
Liouville-type results, see \cite{BM, BuM, CMM, F, FF, M, PQS, RZ,
So, SZ, SZ2} and the references therein; The method of Hardy-Sobolev
inequalities, see \cite{BT, CFM$_2$, C, CFS, GW}. For the detailed
comments of the above methods and the advantages of the bootstrap
methods, we refer to \cite{QS}, see also a survey paper
\cite{S$_2$}.

A similar theorem for system (\ref{sys:main}) with $f,g$ satisfying
(\ref{ass:QS}) was proved by Quittner \& Souplet \cite[Theorem
2.1]{QS}. Based on their a priori estimates, they obtained new
existence theorems for various classes of elliptic systems.

Theorem \ref{thm:a priori estimates} in hand, we are able to obtain
more general existence theorems for system (\ref{sys:main}).
Consider the system (\ref{sys:main}), subject to (\ref{ass:f and g})
and the superlinearity condition
\begin{eqnarray}\label{ass:superlinear}
   f+g\geq \lambda(u+v)-C_1,\ \ \ u, v\geq 0,\
   x\in\Omega,
\end{eqnarray}
where $C_1>0,\ \lambda>\lambda_1$, the first eigenvalue of $-\Delta$
in $H_0^1(\Omega)$.

\begin{thm}\label{thm:main existence} Assume that $f,g$ satisfy $(\ref{ass:f and g})$ and $(\ref{ass:superlinear})$
with $pq>(1-r)(1-s)$ and $(\ref{ass:optimal condition for very
weak})$. Then
\begin{enumerate}
  \item[(a)] any nonnegative $L^1_\delta$-solution $(u,v)$ of
             $(\ref{sys:main})$
             belongs to $L^\infty(\Omega)$ and satisfies
             the a priori estimate $(\ref{est:a priori estimates})$;
  \item[(b)] system $(\ref{sys:main})$
             admits a positive $L^1_\delta$-solution $(u,v)$ if in
             addition $f,g$ satisfy
             \begin{equation}\label{ass:f and g in addition}
                  f+g=o(u+v),\ \ \mathrm{as}\ u,v\rightarrow 0^+,
             \end{equation}
             uniformly in $x\in\Omega$.
\end{enumerate}
\end{thm}

\begin{rem}
If $pq>(1-r)(1-s)$ and $\max\{\alpha,\beta\}>\frac{n-1}{2}$ are
replaced by $pq\leq(1-r)(1-s)$ and
$pq-(1-r)(1-s)<\frac{2}{n-1}\max\{p+1-s,q+1-r\}$ respectively, then
the conclusions of Theorem \ref{thm:a priori estimates} and
\ref{thm:main existence} also hold.
\end{rem}

\begin{rem}
Consider system (\ref{sys:main}) with boundary conditions of the
form $u_\nu=au$, $v_\nu=bv$, where $a, b\in \mathbb{R}$ and $u_\nu$
denotes the derivative of $u$ with respect to the outer unit normal
on $\partial\Omega$. If, for example, $f, g$ satisfy
\begin{eqnarray*}
   f+g\geq C_1(\lambda_1(a)u+\lambda_1(b)v)-C_2,\ \ \ u, v\geq 0,\ x\in\Omega,
\end{eqnarray*}
where $C_1>1$, $C_2\geq 0$ and $\lambda_1(a)$ denotes the first
eigenvalue of $-\Delta$ with boundary conditions $u_\nu=au$, then it
is easy to deduce that
\begin{equation*}
    \|u\|_{L^1}+\|v\|_{L^1}\leq M,
\end{equation*}
with $M$ independent of $u,v$. The proof of Theorem
\ref{thm:abstract} (in Section 2) implies (\ref{est:a priori
estimates}). Using this a priori estimate, we also have a similar
existence theorem of $L^1$-solutions of system (\ref{sys:main}) with
Neumann conditions as Theorem \ref{thm:main existence}.
\end{rem}

If $r=s=0$, under assumptions (\ref{ass:QS}),
(\ref{ass:superlinear}), the system (\ref{sys:main}) was studied by
several authors. Using another bootstrap method, similar results as
the above theorem was obtained in \cite[Theorem 1.1]{QS}, see also
\cite{CFM$_2$, F, FY, Zou} for more related results.

\vskip 3mm

The second existence theorem is about the system
\begin{eqnarray}\label{sys:nuclear reactor}
   &-\Delta u=a(x)u^rv^p-c(x)u,&\ \ {\rm in}\ \Omega, \nonumber\\
   &-\Delta v=b(x)u^qv^s-d(x)v,&\ \ {\rm in}\ \Omega, \\
   &u=v=0,&\ \ {\rm on}\ \partial\Omega,\nonumber
\end{eqnarray}
where $r,s\leq 1$, $pq>(1-r)(1-s)$, $a,b,c,d\in L^\infty(\Omega)$,
$a,b\geq 0$, $\int_\Omega a,\int_\Omega b>0$,
$\inf\{\mathrm{spec}(-\Delta + c)\}>0$, $\inf\{\mathrm{spec}(-\Delta
+ d)\}>0$.

\begin{thm}\label{thm:secondary existence} Assume that
\begin{eqnarray}\label{ass:optimal condition for secondary}
                  \displaystyle\max\{\alpha,\beta\}>\frac{n-1}{2}.
\end{eqnarray}
Then
\begin{enumerate}
  \item[(a)] any nonnegative $L^1_\delta$-solution $(u,v)$ of
             $(\ref{sys:nuclear reactor})$
             belongs to $L^\infty(\Omega)$ and satisfies
             the a priori estimate $(\ref{est:a priori estimates})$;
  \item[(b)] system $(\ref{sys:nuclear reactor})$
             admits a positive $L^1_\delta$-solution $(u,v)$.
\end{enumerate}
\end{thm}

From the above theorem, we obtain the existence theorem for system
(\ref{sys:secondary}).

\begin{cor}
Assume that $r,s\leq 1$, $pq>(1-r)(1-s)$ and $(\ref{ass:optimal
condition for secondary})$ holds. Then system
$(\ref{sys:secondary})$ admits a positive classical solution
$(u,v)$.
\end{cor}

A similar existence result was proved in \cite[Theorem 1.4]{QS} but
under more stronger assumptions. Set
\begin{eqnarray*}
    &\displaystyle\hat{p}=\frac{(n+1)p}{n+1-(n-1)r},\ \ \hat{q}=\frac{(n+1)q}{n+1-(n-1)s},\\
    &\displaystyle\hat{\alpha}=\frac{\hat{p}+1}{(\hat{p}\hat{q}-1)_+},\ \
    \hat{\beta}=\frac{\hat{q}+1}{(\hat{p}\hat{q}-1)_+}.
\end{eqnarray*}
Instead of (\ref{ass:optimal condition for secondary}), they
required that $\max\{\hat{\alpha},\hat{\beta}\}>(n-1)/2$. The a
priori estimates and existence of positive solutions for
(\ref{sys:secondary}) was studied in \cite{CFMT} in the case when
$\Omega=B_R(0)$ and the parameters satisfy $0\leq r,s\leq 1,\ pq\geq
(1-r)(1-s)$, plus some additional conditions. Note that the results
there also cover the case when the Laplace operators are replaced by
$\Delta_m u,\Delta_n u,\ m, n>1$. We refer to \cite{M, RZ, TV, Zh,
B} for existence/nonexistence results for (\ref{sys:secondary}) and
to \cite{DE, Li, W, Zh} and the references therein for related
results on the associated parabolic systems.

\vskip 3mm

\begin{rem}
It was shown in \cite{RZ} that system (\ref{sys:secondary}) has no
positive solutions if $p,q,r,s\geq 1$, $\min\{p+r,q+s\}\geq
(n+2)/(n-2)_+$ and $\Omega$ is star-shaped. It was also proved in
\cite{Zou$_2$} that system (\ref{sys:secondary}) has a positive
solution if $r,s\geq 1$, $pq>(r-1)(s-1)$ and
\begin{equation}\label{ass:Zou}
    \max\{p+r,q+s\}\leq (n+2)/(n-2)_+,
\end{equation}
see also \cite{Zou}. Our result is that system (\ref{sys:secondary})
has a positive solution if $0\leq r,s\leq 1$, $pq>(1-r)(1-s)$ and
(\ref{ass:optimal condition for secondary}) holds. If $r=s=0$, for
the existence of positive solutions of the system
(\ref{sys:secondary}), we have the optimal condition
$\alpha+\beta>(n-2)/2$, see Remark \ref{rem:hyperbola}. We would
like to point out that
\begin{enumerate}
  \item[(i)]  $\max\{p+r,q+s\}\leq (n+1)/(n-1)$ implies
              (\ref{ass:optimal condition for secondary}), but (\ref{ass:Zou})
              does not;
  \item[(ii)] (\ref{ass:optimal condition for secondary}) is much
              more general than (\ref{ass:Zou}). (\ref{ass:optimal condition for secondary})
              allows very large $p$ or $q$;
  \item[(iii)]  If $r=s=0$, (\ref{ass:optimal condition for secondary})
             is stronger than $\alpha+\beta>(n-2)/2$.
\end{enumerate}
So it is still a widely open question what should be the optimal
conditions on $p,q,r,s,n$ for existence of positive solutions to
system (\ref{sys:secondary}).
\end{rem}

\vskip 3mm

A special case of (\ref{sys:nuclear reactor}) is the following model
of a nuclear reactor
\begin{eqnarray}\label{sys:W}
\begin{array}{ll}
  -\Delta u=uv-au, &  {\rm in}\ \Omega, \\
  -\Delta v=bu, &  {\rm in}\ \Omega, \\
  \ \ \ u=v=0, & {\rm on}\ \partial\Omega,
\end{array}
\end{eqnarray}
where $u,v$ present the neutron flux and the temperature,
respectively. This system and the corresponding parabolic system
were studied in \cite{Ch, GW, GW$_2$, Q, QS, QS$_2$}. In
\cite{GW$_2$}, the existence and a priori estimate were obtained
under the assumption $n\leq 3$, or $\Omega$ convex and $n\leq 5$. In
\cite[Theorem 1.2]{QS} and \cite[Theorem 31.17]{QS$_2$}, the
existence and a priori estimate were obtained under weaker
assumption $n\leq 4$ without assuming $\Omega$ convex. Our theorem
recover their result since $\max\{\alpha,\beta\}=2>(n-1)/2$ implies
$n<5$.

\vskip 3mm

In next section, we present our bootstrap procedure. In Section 3,
we prove Theorem \ref{thm:H1 solution}-\ref{thm:very weak solution}.
In Section 4, we prove Theorem \ref{thm:a priori
estimates}-\ref{thm:secondary existence}.

\section{The Bootstrap Procedure}

In what follows we give the definitions of three types of weak
solutions of system (\ref{sys:main}), see \cite[Chapter I]{QS$_2$}.

\begin{defn}
\begin{enumerate}
  \item[(i)] By an $H_0^1$-solution of system (\ref{sys:main}), we mean a couple $(u,v)$ with
\begin{eqnarray*}
    u, v\in H_0^1(\Omega),\ \
    f, g\in H^{-1}(\Omega),
\end{eqnarray*}
satisfying
\begin{eqnarray*}
     &\displaystyle\int_\Omega \nabla u\cdot\nabla\varphi=\int_\Omega f\varphi,\ \
         \int_\Omega \nabla v\cdot\nabla\varphi=\int_\Omega
         g\varphi,\\
     &\mathrm{for\
        all}\ \varphi\in H_0^1(\Omega).
\end{eqnarray*}
  \item[(ii)] By an $L^1$-solution of system (\ref{sys:main}), we mean a couple $(u,v)$ with
\begin{equation*}
    u, v\in L^1(\Omega),\ \ f, g\in
    L^1(\Omega),
\end{equation*}
satisfying
\begin{eqnarray}
     &\displaystyle-\int_\Omega u\Delta\varphi=\int_\Omega f\varphi,\ \
         -\int_\Omega v\Delta\varphi=\int_\Omega
         g\varphi,\nonumber\\
         [-1.5ex]\label{defn:L1 solution}\\[-1.5ex]
     &\mathrm{for\
        all}\ \varphi\in C^2(\overline{\Omega}),\
        \varphi|_{\partial\Omega}=0.\nonumber
\end{eqnarray}
  \item[(iii)] Set $\delta(x):=\mathrm{dist}(x,\partial\Omega)$ and
         $L^1_\delta(\Omega):=L^1(\Omega;\delta(x)\mathrm{d}x)$.
         By an $L^1_\delta$-solution of system (\ref{sys:main}), we mean a couple $(u,v)$ with
\begin{equation*}
    u, v\in L^1(\Omega),\ \ f, g\in
    L^1_\delta(\Omega),
\end{equation*}
satisfying (\ref{defn:L1 solution}).
\end{enumerate}
\end{defn}

The three types of weak solutions of the scalar equation
(\ref{eq:main}) and the linear equation
\begin{eqnarray}\label{eq:linear}
   -\Delta u=\phi,&\ \ {\rm in}\ \Omega;\ \
   u=0,&\ \ {\rm on}\ \partial\Omega,
\end{eqnarray}
are defined similarly. According to \cite[Lemma 1]{BCMR}, if
$\phi\in L^1_\delta(\Omega)$, (\ref{eq:linear}) admits a unique
$L^1_\delta$-solution $u\in L^1(\Omega)$. Moreover, $\|u\|_{L^1}\leq
C\|\phi\|_{L^1_\delta}$ and $\phi\geq 0$ a.e. implies $u\geq 0$ a.e.

The most important regularity results for $L^1$-solutions of the
linear equation (\ref{eq:linear}) is the following
$L^m$-$L^k$-estimates.

\begin{prop}{\rm (see for instance \cite[Proposition 47.5]{QS$_2$})}\label{prop:L1}
Let $1\leq m\leq k\leq \infty$ satisfy
\begin{equation}\label{ineq:m k L1}
   \frac{1}{m}-\frac{1}{k}<\frac{2}{n}.
\end{equation}
Let $u\in L^1(\Omega)$ be the unique $L^1$-solution of
$(\ref{eq:linear})$. If $\phi\in L^m(\Omega)$, then $u\in
L^k(\Omega)$ and satisfies the estimate $\|u\|_{L^k}\leq
C(\Omega,m,k)\|\phi\|_{L^m}$.
\end{prop}

It is well known that the condition (\ref{ineq:m k L1}) is optimal.
For example, let $\Omega=B_1$ be the unit ball. For $1\leq m<k\leq
\infty$, let $n/k<\theta<n/m-2$, which follows from $1/m-1/k>2/n$.
Then $U(r)=r^{-\theta}-1$ is the unique $L^1$-solution of $-\Delta
U=\phi:=\theta(n-\theta-2)r^{-\theta-2}$. But $\phi\in L^m(B_1)$ and
$U\notin L^k(B_1)$, see also \cite[Chapter I]{QS$_2$}.

Obviously, Proposition \ref{prop:L1} holds for the $H_0^1$-solution
of (\ref{eq:linear}). But it is not convenient to derive the optimal
condition for $L^\infty$-regularity of the $H_0^1$-solutions of
system (\ref{sys:main}). For our purpose, we develop an
$L^m$-$L^k$-estimate for the $H_0^1$-solution of (\ref{eq:linear}).
It is an invariant of Proposition \ref{prop:L1}. Let $n\geq 3$, set
$2_*:=2n/(n+2)$. It is the conjugate number of the Sobolev imbedding
exponent, $2n/(n-2)$.

\begin{prop}\label{prop:H1}
Let $1\leq m\leq k\leq \infty$ satisfy
\begin{equation}\label{ineq:m k H1}
   \frac{1}{m}-\frac{1}{k}<\frac{4}{n+2}.
\end{equation}
Let $u\in H_0^1(\Omega)$ be the unique $H_0^1$-solution of
$(\ref{eq:linear})$. If $\phi\in L^{2_*m}(\Omega)$, then $u\in
L^{2_*k}(\Omega)$ and satisfies the estimate $\|u\|_{L^{2_*k}}\leq
C(\Omega,m,k)\|\phi\|_{L^{2_*m}}$.
\end{prop}

The above proposition in hand, the $L^\infty$-regularity of the
$H_0^1$-solutions of (\ref{eq:main}) with $|f|\leq C(1+|u|^p)$ with
$1\leq p<p_S$ follows immediately from a simple bootstrap argument.
It is much simpler than the usual proof, see \cite{BK, St, QS$_2$}.

\vskip 3mm

For all $1\leq k\leq\infty$, define the spaces
$L^k_\delta(\Omega)=L^k(\Omega;\delta(x)\mathrm{d}x)$. For $1\leq
k<\infty$, $L^k_\delta(\Omega)$ is endowed with the norm
\begin{equation*}
    \|u\|_{L^k_\delta}=\left(\int_\Omega|u(x)|^k\delta(x)\mathrm{d}x\right)^{1/k}.
\end{equation*}
Note that $L^\infty_\delta(\Omega)=L^\infty(\Omega;\mathrm{d}x)$,
with the same norm $\|u\|_\infty$. For the $L^1_\delta$-solutions,
we have the following regularity result.

\begin{prop}{\rm(see \cite{FSW}, also \cite{QS, QS$_2$})}\label{prop:L1 delta}
Let $1\leq m\leq k\leq \infty$ satisfy
\begin{equation}\label{ineq:m k L1 delta}
   \frac{1}{m}-\frac{1}{k}<\frac{2}{n+1}.
\end{equation}
Let $u\in L^1(\Omega)$ be the unique $L^1_\delta$-solution of
$(\ref{eq:linear})$. If $\phi\in L^m_\delta(\Omega)$, then $u\in
L^k_\delta(\Omega)$ and satisfies the estimate
$\|u\|_{L^k_\delta}\leq C(\Omega,m,k)\|\phi\|_{L^m_\delta}$.
\end{prop}

The condition (\ref{ineq:m k L1 delta}) is optimal, since for $1\leq
m<k\leq \infty$ and $1/m-1/k>2/(n+1)$, there exists $\phi\in
L^m_\delta(\Omega)$ such that $u\notin L^k_\delta(\Omega)$, where
$u$ is the unique $L^1_\delta$-solution of (\ref{eq:linear}), see
\cite[Theorem 2.1]{S}.

\begin{rem}
According to Proposition \ref{prop:L1}-\ref{prop:L1 delta}, the
assumptions of $h$ in Theorem \ref{thm:H1 solution}-\ref{thm:very
weak solution} are natural.
\end{rem}

\vskip 3mm

In order to give a uniform proof of Theorem \ref{thm:H1
solution}-\ref{thm:very weak solution} (i), we write the three
critical exponents $p_S,\ p_{sg},\ p_{BT}$ as $p_c$.  Denote $B^k$
the spaces $L^{2_*k}(\Omega),\ L^k(\Omega),\ L^k_\delta(\Omega)$,
and $\|\cdot\|_{B^k}$ in $B^k$ the norms $\|\cdot\|_{L^{2_*k}},\
\|\cdot\|_{L^k},\ \|\cdot\|_{L^k_\delta}$. Note that (\ref{ineq:m k
L1})-(\ref{ineq:m k L1 delta}) can be written in one form
\begin{equation}\label{ineq:m k in one form}
   \frac{1}{m}-\frac{1}{k}<\frac{1}{p'_c},
\end{equation}
where $1/p'_c+1/p_c=1$. The optimal conditions  of
$L^\infty$-regularity in Theorem \ref{thm:H1 solution}-\ref{thm:very
weak solution} (i) can also be written in one form
\begin{eqnarray}\label{ass:optimal condition in one form}
                  &\displaystyle\max\{\alpha,\beta\}>\frac{1}{p_c-1},\ \ r,s,\gamma, \sigma<p_c,\nonumber\\
                  [-1.5ex]\\[-1.5ex]
                  &\min\{p+r,q+s\}<p_c,\ h\in B^\theta,\
                  \theta>p'_c.\nonumber
\end{eqnarray}
We shall prove the following theorem.

\begin{thm}\label{thm:abstract}
Assume that $f,g$ satisfy $(\ref{ass:f and g})$ with
$(\ref{ass:optimal condition in one form})$. Then there exists $C>0$
such that for any $(H_0^1,\ L^1,\ L^1_\delta)$-solution $(u,v)$ of
$(\ref{sys:main})$ satisfying
\begin{equation}\label{ass:abstract}
    \|u\|_{B^k}+\|v\|_{B^k}\leq
    M_1(k), \ \ \ \mathrm{for\ all}\ 1\leq k<p_c,
\end{equation}
it follows that $u,v\in L^\infty(\Omega)$ and
\begin{equation*}
    \|u\|_{L^\infty}+\|v\|_{L^\infty}\leq C.
\end{equation*}
The constant $C$ depends only on
$M_1(k),\Omega,p,q,r,s,\gamma,\sigma,C_1$.
\end{thm}

Without loss of generality, we assume that $q+s\geq p+r$. Then
$\beta\geq\alpha$. From (\ref{ass:optimal condition in one form}),
we have
\begin{equation}\label{ass:optimal condition reduced}
        \beta>\frac{1}{p_c-1},
\end{equation}
and
\begin{equation}\label{ass:p+r reduced}
     p+r<p_c.
\end{equation}
\begin{rem}\label{rem:p+r<pc}
If $r\leq 1$, (\ref{ass:p+r reduced}) can be deduced by
(\ref{ass:optimal condition reduced}). In fact, we have
\begin{equation*}
    p+r-1\leq\frac{pq-(1-r)(1-s)}{q+1-r}=\frac{1}{\beta}<p_c-1.
\end{equation*}
\end{rem}

We first prove two lemmas, which assert that by bootstrap only on
the first equation of system (\ref{sys:main}), the integrability of
$u$ can be improved to such an extent that the bootstrap on the
second equation is possible. In the following,
$C=C(M_1,r,s,p,q,\gamma,\sigma,\Omega,C_1)$ is different from line
to line, but it is independent of $(u,v)$ satisfying
(\ref{ass:abstract}). For simplicity, we denote by $|\cdot|_k$ the
norm $\|\cdot\|_{B^k}$.

\begin{lem}\label{lem:p is small}
Let $f,g$ satisfy $(\ref{ass:f and g})$ with $(\ref{ass:optimal
condition in one form})$. If
\begin{equation}\label{ass:p is small}
      p<p_c/p'_c,
\end{equation}
then $|u|_{\infty}\leq C$.
\end{lem}

\begin{proof}
We shall carry out the bootstrap only on the first equation of
system (\ref{sys:main}) to prove $|u|_{\infty}\leq C$.

\textbf{Case I.} $r<1$.

Thanks to (\ref{ass:optimal condition in one form}), (\ref{ass:p+r
reduced}) and (\ref{ass:p is small}) there exists $k$ such that
\begin{equation}\label{res:p is small}
        (p+r)\vee\gamma<k<p_c,\ \ \ \frac{p}{k}<\frac{1}{p'_c}.
\end{equation}
For such $k$ fixed, there exists $\varepsilon>0$ small enough to
satisfy
\begin{equation}\label{btrp:gamma E1}
        \frac{\gamma}{k+m\varepsilon}-\frac{1}{k+(m+1)\varepsilon}<\frac{1}{p'_c},\ \
        \ \mathrm{for\ any\ integer}\ m\geq 0,
\end{equation}
and
\begin{equation}\label{res:r<1}
        r<\frac{k}{k+\varepsilon},
\end{equation}
since $r<1$. From (\ref{res:p is small}) and (\ref{res:r<1}), we
have
\begin{equation}\label{btrp:p r E1}
    \frac{r}{k+m\varepsilon}+\frac{p}{k}-\frac{1}{k+(m+1)\varepsilon}<\frac{1}{p'_c},\ \
        \ \mathrm{for\ any\ integer}\ m\geq 0.
\end{equation}

For $m\geq 0$, set
\begin{equation*}
    \frac{1}{\rho_m}=\frac{r}{k+m\varepsilon}+\frac{p}{k}<1, \ \
    \frac{1}{\varrho_m}=\frac{\gamma}{k+m\varepsilon}<1.
\end{equation*}
From (\ref{res:p is small}), when $m$ is large enough, we have
$\rho_m\wedge\varrho_m>p'_c$. Denote
$m_0=\min\{m:\rho_m\wedge\varrho_m>p'_c\}$. We claim that after
$m_0$-th bootstrap on the first equation, we arrive at the desired
result $|u|_{\infty}\leq C$.

According to (\ref{ass:abstract}), we have $|u|_k\leq C,\ |v|_k\leq
C$. If $m_0=0$, we can take $k$ such that
$p'_c<\rho_0\wedge\varrho_0=k/[(p+r)\vee\gamma]\leq\theta$ and
(\ref{res:p is small}) holds. Then applying Proposition
\ref{prop:L1}-\ref{prop:L1 delta}, using the first equation of
system (\ref{sys:main}), we obtain
\begin{eqnarray}\label{est:epsilon}
     |u|_\infty&\leq& C|f|_{\rho_0\wedge\varrho_0}\nonumber\\
       &\leq& C(||u|^r|v|^p|_{\rho_0\wedge\varrho_0}
               +||u|^\gamma|_{\rho_0\wedge\varrho_0})+|h|_{\rho_0\wedge\varrho_0}\nonumber\\
       &\leq& C(||u|^r|v|^p|_{\rho_0}
               +||u|^\gamma|_{\varrho_0}+1)\nonumber\\
       &\leq& C(|u|_{k}^r|v|_k^p+|u|_{k}^\gamma+1)\nonumber\\
       &\leq& C.
\end{eqnarray}

Now we consider $m_0>0$. If we have got the estimate
$|u|_{k+m\varepsilon}\leq C$ for some $0\leq m<m_0$, then applying
Proposition \ref{prop:L1}-\ref{prop:L1 delta}, using
(\ref{btrp:gamma E1}), (\ref{btrp:p r E1}) and the first equation of
system (\ref{sys:main}), we obtain
\begin{eqnarray}\label{est:m epsilon}
     |u|_{k+(m+1)\varepsilon}&\leq& C|f|_{\rho_m\wedge\varrho_m}\nonumber\\
       &\leq& C(||u|^r|v|^p|_{\rho_m\wedge\varrho_m}
               +||u|^\gamma|_{\rho_m\wedge\varrho_m})+|h|_{\rho_m\wedge\varrho_m}\nonumber\\
       &\leq& C(||u|^r|v|^p|_{\rho_m}
               +||u|^\gamma|_{\varrho_m}+1)\nonumber\\
       &\leq& C(|u|_{k+m\varepsilon}^r|v|_k^p+|u|_{k}^\gamma+1)\nonumber\\
       &\leq& C.
\end{eqnarray}
So we have $|u|_{k+m_0\varepsilon}\leq C$. We can take
$\mathfrak{m}:m_0-1<\mathfrak{m}\leq m_0$ such that
$p'_c<\rho_{\mathfrak{m}}\wedge\varrho_{\mathfrak{m}}\leq \theta$. A
similar argument to (\ref{est:epsilon}) yields $|u|_{\infty}\leq C$.

\textbf{Case II.} $r\geq 1$.

Since $(p+r)\vee\gamma<p_c$, there exist
\begin{eqnarray*}
    &&k:\ (p+r)\vee\gamma<k<p_c,\\
    &&\eta:\eta>1,\ {\rm close\ enough\ to}\ 1,
\end{eqnarray*}
such that
\begin{eqnarray*}
   &\displaystyle\frac{r}{k}+\frac{p}{k}-\frac{1}{\eta
       k}<\frac{1}{p'_c},\\
   &\displaystyle\frac{\gamma}{k}-\frac{1}{\eta
       k}<\frac{1}{p'_c},
\end{eqnarray*}
from which we obtain
\begin{eqnarray*}
   &&\frac{r}{\eta^mk}+\frac{p}{k}-\frac{1}{\eta^{m+1}k}<\frac{2}{n+1},\\
   &&\frac{\gamma}{\eta^mk}-\frac{1}{\eta^{m+1}k}<\frac{2}{n+1},
\end{eqnarray*}
for any integer $m\geq 0$. Similarly to the arguments of Case I, we
also have $|u|_{\infty}\leq C$.

The proof of the lemma is complete.
\end{proof}

\begin{lem}\label{lem:p is large}
Let $f,g$ satisfy $(\ref{ass:f and g})$ with $(\ref{ass:optimal
condition in one form})$. If
\begin{equation}\label{ass:p is large}
      p\geq p_c/p'_c.
\end{equation}
Let $k^*: p_c<k^*\leq\infty$ be the solution of
\begin{equation}\label{def:k *}
       \frac{r}{k^*}+\frac{p}{p_c}-\frac{1}{k^*}=\frac{1}{p'_c}.
\end{equation}
Then for any $1\leq k_1<k^*$, we have $|u|_{k_1}\leq C$.
\end{lem}

\begin{proof}
According to (\ref{ass:p+r reduced}) and (\ref{ass:p is large}), we
necessarily have $r<1$. We shall also carry out the bootstrap only
on the first equation of system (\ref{sys:main}) to prove
$|u|_{k_1}\leq C$. We first consider the case where $p>p_c/p'_c$. So
$p_c<k^*<\infty$. For any $\varepsilon: 0<\varepsilon\ll 1$, set
$k_{\varepsilon}=k^*-\varepsilon$. Thanks to (\ref{ass:p+r reduced})
and (\ref{def:k
*}), since $r<1$, there exist
\begin{eqnarray*}
   &&k:\ (p+r)\vee\gamma<k<p_c, \ \ \mathrm{close\ enough\ to}\ p_c,\\
   &&\tau:\ r<\tau<1,\ {\rm close\ to}\ 1,
\end{eqnarray*}
such that
\begin{eqnarray}
    &&\displaystyle\frac{r}{k_{\varepsilon}}+\frac{p}{k}-\frac{1}{k_{\varepsilon}}<\frac{1}{p'_c},
      \label{res:k *-varepsilon}\\
    &&rk_{\varepsilon}^\tau<\tau k,\label{res:exist tau}\\
    &&\frac{\gamma}{k}-
          \frac{1}{k_{\varepsilon}^\tau}<\frac{1}{p'_c},\label{res:1 tau gamma}
\end{eqnarray}
where
$k_{\varepsilon}^{\tau^m}=k_{\varepsilon}-\tau^m(k_{\varepsilon}-k)$
for $m\geq 0$. In fact, (\ref{res:k *-varepsilon}) is a small
perturbation of (\ref{def:k *}) with respect to $k^*$ and,
(\ref{res:exist tau}) is a small perturbation of itself with
$\tau=1$. A careful computation yields that
\begin{eqnarray*}
     &&\frac{r}{k_{\varepsilon}^{\tau^m}}-\frac{1}{k_{\varepsilon}^{\tau^{m+1}}}
       <\frac{r}{k_{\varepsilon}}-\frac{1}{k_{\varepsilon}},\ \
           \ \mathrm{for\ any\ integer}\ m\geq 0,\ \ \ \ (\mathrm{using}\ (\ref{res:exist tau}))\\
     &&\frac{\gamma}{k_{\varepsilon}^{\tau^m}}-\frac{1}{k_{\varepsilon}^{\tau^{m+1}}}
       <\frac{\gamma}{k_{\varepsilon}^{\tau^{m-1}}}-\frac{1}{k_{\varepsilon}^{\tau^m}},\ \
           \ \mathrm{for\ any\ integer}\ m\geq 1.\ \ \ \ (\mathrm{using}\
           \gamma\geq 1)
\end{eqnarray*}
So, according to (\ref{res:k *-varepsilon}) and (\ref{res:1 tau
gamma}), we have
\begin{eqnarray}
        &&\frac{r}{k_{\varepsilon}^{\tau^m}}+\frac{p}{k}
           -\frac{1}{k_{\varepsilon}^{\tau^{m+1}}}<\frac{1}{p'_c},\label{bt:tau m}\\
        &&\frac{\gamma}{k_{\varepsilon}^{\tau^m}}
           -\frac{1}{k_{\varepsilon}^{\tau^{m+1}}}<\frac{1}{p'_c}.\label{bt:tau m gamma}
\end{eqnarray}
for any integer $m\geq 0$.

Set
\begin{equation*}
    \frac{1}{\rho_m}=\frac{r}{k_{\varepsilon}^{\tau^m}}+\frac{p}{k}<1, \ \
    \frac{1}{\varrho_m}=\frac{\gamma}{k_{\varepsilon}^{\tau^m}}<1.
\end{equation*}
Note that
\begin{eqnarray*}
    \frac{1}{\vartheta}=\frac{r}{k_{\varepsilon}}+\frac{p}{k}>\frac{r}{k^*}+\frac{p}{p_c}\geq\frac{1}{p'_c}.
\end{eqnarray*}
So $\rho_m\wedge\varrho_m<\vartheta<p'_c<\theta$. Then
$|h|_{\rho_m\wedge\varrho_m}\leq C|h|_\theta\leq C$ for all $m\geq
0$.

We already have $|u|_k\leq C,\ |v|_k\leq C$ from
(\ref{ass:abstract}). If we have got
$|u|_{k_{\varepsilon}^{\tau^m}}\leq C$ for some $m\geq 0$, applying
Proposition \ref{prop:L1}-\ref{prop:L1 delta}, using (\ref{bt:tau
m}), (\ref{bt:tau m gamma}) and the first equation of system
(\ref{sys:main}), similarly to (\ref{est:m epsilon}), we obtain
$|u|_{k_{\varepsilon}^{\tau^{m+1}}}\leq C$. So, for any integer
$m\geq 0$, there holds $|u|_{k_{\varepsilon}^{\tau^m}}\leq C$.
Noting that $k_{\varepsilon}^{\tau^m}\rightarrow k_{\varepsilon}$ as
$m\rightarrow\infty$, we prove the lemma for $p>p_c/p'_c$.

If $p=p_c/p'_c$, we have $k^*=\infty$. The above proof is also valid
when $k_\varepsilon$ is replaced by any arbitrary large number. The
proof is complete.
\end{proof}

Lemma \ref{lem:p is small} and \ref{lem:p is large} in hand, we can
prove Theorem \ref{thm:abstract}.

\vskip 3mm

\noindent\textbf{Proof of Theorem \ref{thm:abstract}.}

\textbf{Case I.} $p<p_c/p'_c$.

According to Lemma \ref{lem:p is small}, $|u|_{\infty}\leq C$. Since
$s,\sigma<p_c$, a simple bootstrap argument on the second equation
yields that $|v|_{\infty}\leq C$.

\textbf{Case II.} $p=p_c/p'_c$.

According to Lemma \ref{lem:p is large}, $|u|_{k_1}\leq C$ for any
$k_1\geq 1$. Take $k_1$ large enough and $k:s\vee\sigma<k<p_c$ such
that
\begin{equation*}
  \frac{q}{k_1}<\frac{1}{p'_c},\ \ \frac{q}{k_1}+\frac{s}{k}<1.
\end{equation*}
Similarly to the proof of Lemma \ref{lem:p is small}, we have
$|v|_{\infty}\leq C$. So we also have $|u|_{\infty}\leq C$ since
$r,\gamma<p_c$.

\textbf{Case III.} $p>p_c/p'_c$. In this case we necessarily have
$r<1$.

According to (\ref{ass:optimal condition reduced}) and (\ref{def:k
*}), there exist
\begin{eqnarray*}
    &&k_1:\ p_c<k_1<k^*,\ {\rm close\ enough\ to}\ k^*,\\
    &&k:\ (p+r)\vee\gamma\vee\sigma<k<p_c,\ \ \mathrm{close\ enough\ to}\ p_c,\\
    &&\eta:\ \eta>1\ {\rm close\ enough\ to}\ 1,
\end{eqnarray*}
such that
\begin{eqnarray}
   &&\displaystyle\frac{q}{k_1}+\frac{s}{k}<1,\label{bt:condition for bootstrap}\\
   &&\displaystyle\frac{r}{k_1}+\frac{p}{\eta k}-\frac{1}{\eta k_1}<\frac{1}{p'_c},\label{bt:first step 1}\\
   &&\displaystyle\frac{q}{k_1}+\frac{s}{k}-\frac{1}{\eta k}<\frac{1}{p'_c},\label{bt:first step 2}\\
   &&\displaystyle\frac{\gamma}{k_1}-\frac{1}{\eta k_1}<\frac{1}{p'_c},\label{bt:bootstrap for gamma 1}\\
   &&\displaystyle\frac{\sigma}{k}-\frac{1}{\eta k}<\frac{1}{p'_c}.\label{bt:bootstrap for sigma 1}
\end{eqnarray}
In fact, (\ref{bt:first step 2}) is equivalent to (\ref{bt:condition
for bootstrap}). (\ref{bt:condition for bootstrap}) with $k_1=k^*$
and $k=p_{BT}$ is exactly (\ref{ass:optimal condition reduced}). So,
(\ref{bt:condition for bootstrap}) and (\ref{bt:first step 2}) are
just small perturbations of (\ref{ass:optimal condition reduced}).
(\ref{bt:first step 1}) is a small perturbation of (\ref{def:k *}).
Multiplying the LHS of (\ref{bt:first step 1})-(\ref{bt:bootstrap
for sigma 1}) by $1/\eta^m$, we have
\begin{eqnarray}
   &&\displaystyle\frac{r}{\eta^m k_1}+\frac{p}{\eta^{m+1} k}-
       \frac{1}{\eta^{m+1} k_1}<\frac{1}{p'_c},\ \
       \frac{\gamma}{\eta^m k_1}-\frac{1}{\eta^{m+1} k_1}<\frac{1}{p'_c},\label{bt:m step 1}\\
   &&\displaystyle\frac{q}{\eta^m k_1}+\frac{s}{\eta^m k}
       -\frac{1}{\eta^{m+1} k}<\frac{1}{p'_c},\ \
       \frac{\sigma}{\eta^m k}-\frac{1}{\eta^{m+1} k}<\frac{1}{p'_c},\label{bt:m step 2}
\end{eqnarray}
for any integer $m\geq 0$.

Set
\begin{eqnarray*}
  &&\frac{1}{\mu_m}=\frac{r}{\eta^m k_1}+\frac{p}{\eta^{m+1} k}<1, \ \
    \frac{1}{\nu_m}=\frac{\gamma}{\eta^m k_1}<1,\\
  &&\frac{1}{\rho_m}=\frac{q}{\eta^m k_1}+\frac{s}{\eta^m k}<1, \ \
    \frac{1}{\varrho_m}=\frac{\sigma}{\eta^m k}<1.
\end{eqnarray*}
Since $\eta>1$, for $m$ large enough, we have
$\rho_m\wedge\varrho_m>p'_c$ and $\mu_m\wedge\nu_m>p'_c$. Denote
$m_0=\min\{m:(\rho_m\wedge\varrho_m)\vee(\mu_m\wedge\nu_m)>p'_c\}$.
We may assume that $\rho_{m_0}\wedge\varrho_{m_0}>p'_c$. We claim
that after $m_0$-th alternate bootstrap on system (\ref{sys:main}),
we shall arrive at the desired result $|v|_{\infty}\leq C$.

We already have $|u|_{k_1}\leq C$ (from Lemma \ref{lem:p is large})
and $|v|_k\leq C$ (from (\ref{ass:abstract})). If $m_0=0$, we can
take $k,\ k_1$ such that $p'_c<\rho_0\wedge\varrho_0\leq\theta$.
Then applying Proposition \ref{prop:L1}-\ref{prop:L1 delta}, using
the second equation of system (\ref{sys:main}), a similar argument
to (\ref{est:epsilon}) yields that $|v|_{\infty}\leq C$. So we also
have $|u|_{\infty}\leq C$ since $r,\gamma<p_c$.

Now we consider $m_0>0$. If we have got the estimate $|u|_{\eta^m
k_1}+|v|_{\eta^m k}\leq C$ for some $0\leq m<m_0$, then applying
Proposition \ref{prop:L1}-\ref{prop:L1 delta}, using (\ref{bt:m step
2}) and the second equation of system (\ref{sys:main}), a similar
argument to (\ref{est:m epsilon}) yields that $|v|_{\eta^{m+1}k}\leq
C$. Then using (\ref{bt:m step 1}) and the first equation of system
(\ref{sys:main}), we obtain $|u|_{\eta^{m+1}k}\leq C$. So we have
$|u|_{\eta^{m_0}k}+|v|_{\eta^{m_0}k}\leq C$. We can take
$\mathfrak{m}:m_0-1<\mathfrak{m}\leq m_0$ such that
$p'_c<\rho_{\mathfrak{m}}\wedge\varrho_{\mathfrak{m}}\leq \theta$. A
similar argument to (\ref{est:epsilon}) yields $|v|_{\infty}\leq C$.
So we also have $|u|_{\infty}\leq C$ since $r,\gamma<p_c$. The proof
is complete.\hskip 50mm  $\Box$

Theorem \ref{thm:abstract} also holds if $pq\leq (1-r)(1-s)$ in
(\ref{ass:f and g}), we have the following theorem.

\begin{thm}\label{thm:abstract pq leq (1-r)(1-s)}
Assume that $f,g$ satisfy $(\ref{ass:f and g})$ with $pq\leq
(1-r)(1-s)$ and $(\ref{ass:optimal condition in one form})$ where
$\max\{\alpha,\beta\}>1/(p_c-1)$ is replaced by
$pq-(1-r)(1-s)<(p_c-1)\max\{p+1-s,q+1-r\}$. Then the conclusion of
Theorem \ref{thm:abstract} holds.
\end{thm}

\begin{proof}
Assume that $q+s\geq q+r$. Note that
\begin{equation*}
     \frac{q}{k^*}+\frac{s}{p_c}<1
\end{equation*}
is equivalent to $pq-(1-r)(1-s)<(p_c-1)\max\{p+1-s,q+1-r\}$. So the
proof is essentially word by word the same as the proof of Theorem
\ref{thm:abstract}.
\end{proof}

\section{$L^\infty$-regularity}

In this section, we prove Theorem \ref{thm:H1
solution}-\ref{thm:very weak solution}.

\vskip 2mm

\textbf{Proof of Theorem \ref{thm:H1 solution}.}

(i) If $n=1,2$, the $L^\infty$-regularity of $H_0^1$-solutions
follows directly from the Sobolev imbedding theorem and Proposition
\ref{prop:L1}. If $n\geq 3$, since $u,v\in H_0^1(\Omega)$, we have
(\ref{ass:abstract}) from the Sobolev imbedding theorem. Then the
$L^\infty$-regularity follows from Theorem \ref{thm:abstract} with
$p_c=(n+2)/(n-2)$ and $B^1=L^{2_*}(\Omega)$ according to
(\ref{ass:optimal condition for H1}).

(ii) Let $(u,v)=(c_1|x|^{-2\alpha}-c_1,c_2|x|^{-2\beta}-c_2)$, where
$c_1,c_2$ are determined by $c_1^{r-1}c_2^p=2\alpha(n-2-2\alpha)$,
$c_1^qc_2^{s-1}=2\beta(n-2-2\beta)$.  Since
$\alpha,\beta<(n-2)/4<(n-2)/2$,  we have $c_1,c_2>0$. Obviously,
\begin{eqnarray*}
    &&-\Delta u=2c_1\alpha(n-2-2\alpha)|x|^{-2\alpha-2}
        =c_1^rc_2^p|x|^{-2\alpha r-2\beta p}=(u+c_1)^r(v+c_2)^p,\\
    &&-\Delta v=2c_2\beta(n-2-2\beta)|x|^{-2\beta-2}
        =c_1^qc_2^s|x|^{-2\alpha q-2\beta s}=(u+c_1)^q(v+c_2)^s.
\end{eqnarray*}
It is easy to verify that $(u,v)$ is an $H_0^1$-solution of system
(\ref{sys:main}) in $B_1$ with
$f=(u+c_1)^r(v+c_2)^p,g=(u+c_1)^q(v+c_2)^s$.\ \ \ \ $\Box$

\vskip 2mm

\textbf{Proof of Theorem \ref{thm:L1 solution}.}

(i) If $n=1,2$, the $L^\infty$-regularity of $L^1$-solutions follows
directly from Proposition \ref{prop:L1}. If $n\geq 3$, since
$f(\cdot,u,v),g(\cdot,u,v)\in L^1(\Omega)$, we have
(\ref{ass:abstract}) from Proposition \ref{prop:L1}. Then the
$L^\infty$-regularity follows from Theorem \ref{thm:abstract} with
$p_c=n/(n-2)$ and $B^1=L^1(\Omega)$ according to (\ref{ass:optimal
condition for L1}).

(ii) Since $\alpha,\beta<(n-2)/2$, $(u,v)$ constructed in the proof
of Theorem \ref{thm:H1 solution} (ii) is also a $L^1$-solution of
system (\ref{sys:main}) in $B_1$ with
$f=(u+c_1)^r(v+c_2)^p,g=(u+c_1)^q(v+c_2)^s$.\ \ \ \ $\Box$

\vskip 2mm

\textbf{Proof of Theorem \ref{thm:very weak solution}.}

(i) If $n=1$, the $L^\infty$-regularity of $L^1_\delta$-solutions
follows directly from Proposition \ref{prop:L1 delta}. If $n\geq 2$,
we have (\ref{ass:abstract}) since $f(\cdot,u,v),g(\cdot,u,v)\in
L^1_\delta(\Omega)$ from Proposition \ref{prop:L1 delta}. Then the
$L^\infty$-regularity follows from Theorem \ref{thm:abstract} with
$p_c=(n+1)/(n-1)$ and $B^1=L^1_\delta(\Omega)$ according to
(\ref{ass:optimal condition for very weak}).

(ii) Assume that $0\in \partial\Omega$. Let $-1<\theta<(n-1)/2$. Let
$\Sigma_1$ be a revolution cone of vertex zero and
$\Sigma:=\Sigma_1\cap B_R\in \Omega$ for sufficiently small $R>0$.
Then $\phi=|x|^{-2(\theta+1)}\mathbf{1}_\Sigma\in
L_\delta^1(\Omega)$ and according to \cite[Lemma 5.1]{S}, the
solution $U>0$ of (\ref{eq:linear}) satisfies $U\geq
C|x|^{-2\theta}\mathbf{1}_\Sigma$. Set
$\phi=|x|^{-2(\alpha+1)}\mathbf{1}_\Sigma$ and
$\psi=|x|^{-2(\beta+1)}\mathbf{1}_\Sigma$, and $u,v>0$ be the
corresponding solutions of (\ref{eq:linear}). We have $u,v\notin
L^\infty$, and
\begin{eqnarray*}
    &&u^rv^p\geq C|x|^{-2\alpha r-2\beta p}\mathbf{1}_\Sigma=C
        |x|^{-2(\alpha+1)}\mathbf{1}_\Sigma=C\phi,\\
    &&u^qv^s\geq C'|x|^{-2\alpha q-2\beta s}\mathbf{1}_\Sigma=C'
        |x|^{-2(\beta+1)}\mathbf{1}_\Sigma=C'\psi.
\end{eqnarray*}
Setting $a(x)=\phi/(u^rv^p)\geq 0$, $b(x)=\psi/(u^qv^s)\geq 0$, we
get
\begin{eqnarray*}
   &-\Delta u=\phi=a(x)u^rv^p,&\ \ {\rm in}\ \Omega, \nonumber\\
   &-\Delta v=\psi=b(x)u^qv^s,&\ \ {\rm in}\ \Omega,
\end{eqnarray*}
and $a(x)\leq 1/C$, $b(x)\leq 1/C'$, hence $a,b\in L^\infty$.\ \ \
$\Box$

\textbf{Proof of Theorem \ref{thm:pq leq (1-r)(1-s)}.}

The proof is word by word the same as the proof of Theorem
\ref{thm:H1 solution}-\ref{thm:very weak solution} (i). \ \ \ $\Box$

\section{A priori estimates of $L^1_\delta$-solutions and existence theorems}

In order to prove Theorem \ref{thm:a priori estimates}, we recall a
special property of the $L^1_\delta$-solutions, which is a
consequence of Proposition \ref{prop:L1 delta}, see
\cite[Proposition 2.2, 2.3]{QS}.

\begin{prop}\label{porp:L k delta}
Let $(u,v)$ be the $L^1_\delta$-solution of system
$(\ref{sys:main})$ with $f,g$ satisfying (\ref{ass:a priori
estimates}) and let $1\leq k<p_{BT}$. Then $u,v\in
L^k_\delta(\Omega)$ and satisfies the estimate
$\|u\|_{L^k_\delta}+\|v\|_{L^k_\delta}\leq
C(\Omega,k,C_2)(\|u\|_{L^1_\delta}+\|v\|_{L^1_\delta}+\|h_1\|_{L^1_\delta})$.
\end{prop}

\begin{proof}
The proof is similar to that of \cite[Proposition 2.2]{QS}. Let
$\varphi_1(x)$ be the first eigenfunction of $-\Delta$ in
$H_0^1(\Omega)$. Recall that
\begin{equation*}
    c_1\delta(x)\leq\varphi_1(x)\leq c_2\delta(x),\ \ x\in \Omega,
\end{equation*}
for some $c_1,c_2>0$. We have
\begin{eqnarray*}
  \int_\Omega(|f|+|g|)\varphi_1&=& \int_\Omega(|\Delta u|+|\Delta v|)\varphi_1=2\int_\Omega((\Delta
              u)_++(\Delta v)_+)\varphi_1-\int_\Omega\varphi_1(\Delta u+\Delta v)\\
              &\leq& 2\int_\Omega(C_2(u_++v_+)+h_+)\varphi_1+\lambda_1\int_\Omega
              (u+v)\varphi_1\\
              &\leq&
              C(\Omega,C_2)(\|u_+\|_{L^1_\delta}+\|v_+\|_{L^1_\delta}+\|h_+\|_{L^1_\delta})\\
              &\leq&
              C(\Omega,C_2)(\|u\|_{L^1_\delta}+\|v\|_{L^1_\delta}+\|h\|_{L^1_\delta}).
\end{eqnarray*}
Applying Proposition \ref{prop:L1 delta} with $m=1$, we have
\begin{equation*}
  \|u\|_{L^k_\delta}+\|v\|_{L^k_\delta}\leq
    C(\Omega,k,C_2)(\|u\|_{L^1_\delta}+\|v\|_{L^1_\delta}+\|h_1\|_{L^1_\delta}).
\end{equation*}
\end{proof}

\textbf{Proof of Theorem \ref{thm:a priori estimates}.}

Since $f,g$ satisfy (\ref{ass:a priori estimates}), from Proposition
\ref{porp:L k delta}, (\ref{ass:abstract}) can be deduced by
(\ref{ass:u v}). So this theorem follows immediately from Theorem
\ref{thm:abstract} with $p_c=(n+1)/(n-1)$ and
$B^1=L^1_\delta(\Omega)$. \ \ \ \ \ $\Box$

\vskip 3mm

From Theorem \ref{thm:a priori estimates}, in order to obtain the a
priori estimate (\ref{est:a priori estimates}), we only have to
obtain, for all $L^1_\delta$-solutions $(u,v)$ of system
(\ref{sys:main}), $\|u\|_{L^1_\delta}+\|v\|_{L^1_\delta}\leq M$ for
some $M$ independent of $u,v$. In the following we give some
propositions which assert the a priori estimate (\ref{est:a priori
estimates}).

\begin{prop}\cite[Proposition 3.1]{QS}\label{prop:QS1}
If $f,g$ satisfy $(\ref{ass:superlinear})$ with $\lambda>\lambda_1$,
then any nonnegative $L^1_\delta$-solution of system
$(\ref{sys:main})$ satisfies $(\ref{ass:u v})$ with $M$ independent
of $u,v$.
\end{prop}

\begin{prop}\cite[Proposition 3.2]{QS}\label{prop:QS2}
If $f,g$ satisfy
\begin{eqnarray}\label{ass:QS2}
   &&f\geq C_1u^rv^p-C_2u,\nonumber\\
        [-1.5ex]&&\hskip 60mm u, v\geq 0,\ x\in\Omega\\[-1.5ex]
   &&g\geq C_1u^qv^s-C_2v,\nonumber
\end{eqnarray}
where $r,s<1,\ pq>(1-r)(1-s)$. Then any nonnegative
$L^1_\delta$-solution of system $(\ref{sys:main})$ in $H_0^1\cap
L^\infty$ satisfies $(\ref{ass:u v})$ with $M$ independent of $u,v$.
\end{prop}

Proposition \ref{prop:QS2} can be extended to some case where
$r,s\geq 1$, see \cite[Proposition 3.5]{QS}, see also \cite[Theorem
1.4 (ii), (iii)]{QS} for the precise assumptions.

The following proposition gives the uniform $L^1_\delta$-estimates
of the $L^1_\delta$-solutions of system (\ref{sys:nuclear reactor})
where $r,s\leq 1$.

\begin{prop}\label{prop:nuclear reactor}
Any nonnegative $L^1_\delta$-solution $(u,v)$ of system
$(\ref{sys:nuclear reactor})$ satisfies $(\ref{ass:u v})$ with $M$
independent of $u,v$.
\end{prop}

\begin{proof}
We use the idea of \cite[Proposition 4.1]{S}. Denote $G(x,y),\
V(x,y)$ the Green functions in $\Omega$ for $-\Delta$ and
$-\Delta+q(x)$. If $\inf\{\mathrm{spec}(-\Delta + q)\}>0$, by
\cite[Theorem 8]{Zhao}, there exists a positive constant
$C=C(\Omega,q)$ such that
\begin{equation*}
     \frac{1}{C}G(x,y)\leq V(x,y)\leq CG(x,y).
\end{equation*}
By \cite[Lemma 3.2]{BC}, we know that
\begin{equation*}
    G(x,y)\geq C\delta(x)\delta(y)\ \ \ \mathrm{for}\
    x,y\in\overline{\Omega}.
\end{equation*}
So we also have
\begin{equation*}
    V(x,y)\geq C\delta(x)\delta(y)\ \ \ \mathrm{for}\
    x,y\in\overline{\Omega},
\end{equation*}
for some constant $C>0$. Denote $\varphi_q(x)$ the first
eigenfunction of $-\Delta+q(x)$ in $H_0^1(\Omega)$ and $\lambda_q$
the first eigenvalue. Recall that
\begin{equation*}
    c_1\delta(x)\leq\varphi_q(x)\leq c_2\delta(x),\ \ x\in \Omega,
\end{equation*}
for some $c_1,c_2>0$. Let $w$ be the solution of the linear equation
\begin{equation*}
     -\Delta w+q(x)w=\phi(x), \ x\in\Omega;\ \ w=0,\
     x\in\partial\Omega.
\end{equation*}
If $\phi\in L_\delta^1$ is nonnegative, then we have
\begin{eqnarray*}
    w=\int_\Omega V(x,y)\phi(x)\geq C(\int_\Omega\phi\delta)\delta\geq
       C(\int_\Omega\phi\varphi_q)\varphi_q
\end{eqnarray*}
with $C$ depending only on $\Omega,q(x)$. Let $(u,v)$ be a
nonnegative $L^1_\delta$-solution of (\ref{sys:nuclear reactor}).
Set
\begin{eqnarray*}
    A=\int_\Omega a(x)u^rv^p\varphi_c,\ \ B=\int_\Omega b(x)u^qv^s\varphi_d.
\end{eqnarray*}
Then we have
\begin{eqnarray*}
    u\geq CA\varphi_c,\ \ v\geq CB\varphi_d.
\end{eqnarray*}
Therefore we obtain
\begin{eqnarray}
  &&A\geq C\int_\Omega
     a\varphi_c^{r+1}\varphi_d^qA^rB^p\geq CA^rB^p,\label{res:A}\\
  &&B\geq C\int_\Omega
     b\varphi_c^q\varphi_d^{s+1}A^qB^s\geq CA^qB^s\label{res:B}.
\end{eqnarray}
If $r=1$ or $s=1$, $A,B\leq C$ obviously. We consider $r<1$. From
(\ref{res:A}), we have $A^{1-r}\geq CB^p$. So combining with
(\ref{res:B}), we obtain $B\geq CB^{pq/(1-r)+s}$. Since
$pq>(1-r)(1+s)$, we have $B\leq C$. From (\ref{res:B}), we also have
$A\leq C$. Using $\varphi_c$ as a testing function in the first
equation of (\ref{sys:nuclear reactor}) and $\varphi_d$ in the
second equation, this yields that
\begin{eqnarray*}
  \int_\Omega u\varphi_c=\int_\Omega a(x)u^rv^p\varphi_c=A\leq C,\\
  \int_\Omega v\varphi_d=\int_\Omega b(x)u^qv^s\varphi_d=B\leq C,.
\end{eqnarray*}
The proof is complete.
\end{proof}

Now we can prove our existence theorems. The proof is standard, see
\cite{QS}. For the readers' convenience, we give the details.

\vskip 3mm

\textbf{Proof of Theorem \ref{thm:main existence}.}

(a) This is a direct consequence of Theorem \ref{thm:a priori
estimates} and Proposition \ref{prop:QS1}.

(b) Let $K$ be the positive cone in $X:=L^\infty(\Omega)\times
L^\infty(\Omega)$ and let $S:X\rightarrow X:(\phi,\psi)\mapsto
(u,v)$ be the solution operator of the linear problem
\begin{eqnarray*}
   &-\Delta u=\phi,\ \ -\Delta v=\psi,\ \ {\rm in}\ \Omega,\\
   &u=v=0,\ \ {\rm on}\ \partial\Omega.
\end{eqnarray*}
Since any nonnegative $L^1_\delta$-solution of (\ref{sys:main}) is
in $L^\infty$ by part (a), the system (\ref{sys:main}) is equivalent
to the equation $(u,v)=T(u,v)$, where $T:X\rightarrow X$ is a
compact operator defined by $T(u,v)=S(f(\cdot,u,v),g(\cdot,u,v))$.
Let $W\subset K$ be relatively open, $Tz\neq z$ for
$z\in\overline{W}\setminus W$, and let $i_K(T,W)$ be the fixed point
index of $T$ with respect to $W$ and $K$ (see \cite{AF} the
definition and basic properties of this index).

If $W_\varepsilon=\{(u,v)\in K:\|(u,v)\|_X<\varepsilon\}$ and
$\varepsilon>0$ is small enough, then (\ref{ass:f and g in
addition}) guarantees $H_1(\mu,u,v)\neq(u,v)$ for any $\mu\in[0,1]$
and $(u,v)\in \overline{W_\varepsilon}\setminus W_\varepsilon$,
where
\begin{equation*}
     H_1(\mu,u,v)=\mu T(u,v)=S(\mu f(\cdot,u,v),\mu g(\cdot,u,v)).
\end{equation*}
Consequently,
\begin{equation*}
    i_K(T,W_\varepsilon)=i_K(H_1(1,\cdot,\cdot),W_\varepsilon)=
      i_K(H_1(0,\cdot,\cdot),W_\varepsilon)=i_K(0,W_\varepsilon)=1.
\end{equation*}

On the other hand, if $R>0$ is large, then our a priori esstimates
guarantee $H_2(\mu,u,v)\neq(u,v)$ for any $\mu\in [0,C_1+1]$ and
$(u,v)\in \overline{W_R}\setminus W_R$, where
\begin{equation*}
     H_2(\mu,u,v)=S(f(\cdot,u,v)+\mu,g(\cdot,u,v)).
\end{equation*}
Using $\varphi_1$ as a testing function we easily see that
$H_2(C_1+1,u,v)=(u,v)$ does not possess nonnegative solutions, hence
\begin{equation*}
    i_K(T,W_R)=i_K(H_2(C_1+1,\cdot,\cdot),W_R)=0.
\end{equation*}
Consequently, $i_K(T,W_R\setminus\overline{W_\varepsilon})=-1$,
which implies existence of a positive solution of (\ref{sys:main}).
The proof is complete.\ \ \ \ $\Box$

\vskip 3mm

\textbf{Proof of Theorem \ref{thm:secondary existence}.}

(a) This is a direct consequence of Theorem \ref{thm:a priori
estimates} and Proposition \ref{prop:nuclear reactor}.

(b) Let $K,X,W_\varepsilon$ be the same as in the proof of Theorem
\ref{thm:main existence} (b), let $S$ be the solution operator of
the linear problem
\begin{eqnarray*}
   &-\Delta u+c(x)u=\phi,\ \ {\rm in}\ \Omega,\\
   &-\Delta v+d(x)v=\psi,\ \ {\rm in}\ \Omega,\\
   &u=v=0,\ \ {\rm on}\ \partial\Omega.
\end{eqnarray*}
Let us show that $H_1(\mu,u,v)\neq(u,v)$ for any $\mu\in[0,1]$ and
$(u,v)\in \overline{W_\varepsilon}\setminus W_\varepsilon$ for
$\varepsilon$ small. Assume by contrary $(u,v)\in
\overline{W_\varepsilon}\setminus W_\varepsilon$,
$H_1(\mu,u,v)=(u,v)$. Then $u\neq 0$, $v\neq 0$ and the standard
$L^z$-estimates (with $z>n/2$) guarantee
\begin{equation*}
     \|u\|_\infty\leq C\|u\|^r_\infty\|v\|^p_\infty,\ \ \|v\|_\infty\leq
     C\|u\|^q_\infty\|v\|^s_\infty.
\end{equation*}
Hence
\begin{equation*}
     \|u\|^{(1-r)(1-s)}_\infty\leq C\|u\|_\infty^{pq},
\end{equation*}
which contradicts $pq>(1-r)(1-s)$ if $\varepsilon$ is small enough.

On the other hand, if $R>0$ is large, then our a priori esstimates
guarantee $H_2(\mu,u,v)\neq(u,v)$ for any $\mu\in [0,\lambda_c]$ and
$(u,v)\in \overline{W_R}\setminus W_R$, where
\begin{equation*}
     H_2(\mu,u,v)=S(f(\cdot,u,v)+\mu(u+1),g(\cdot,u,v)).
\end{equation*}
and $\lambda_c$ is the first eigenvalue of $-\Delta+c(x)$ in
$H_0^1(\Omega)$. Using $\varphi_c$ as a testing function we easily
see that $H_2(\lambda_c,u,v)=(u,v)$ does not possess nonnegative
solutions, hence
\begin{equation*}
    i_K(T,W_R)=i_K(H_2(\lambda_c,\cdot,\cdot),W_R)=0.
\end{equation*}
Consequently, $i_K(T,W_R\setminus\overline{W_\varepsilon})=-1$,
which implies existence of a positive solution of (\ref{sys:nuclear
reactor}). The proof is complete.\ \ \ \ $\Box$

\vskip 3mm

\noindent\emph{Acknowledgements.} Li Yuxiang is grateful to
Professor Philippe Souplet for many helpful discussions and remarks
during the preparation of this paper and, for his warm reception and
many helps when Li visited the second address.



\begin{thebibliography}{999999}
\bibitem[A]{A} Aviles, P.,
            \textit{On isolated singularities in some nonlinear partial differential equations},
             Indiana Univ. Math. J.
             \textbf{32} (1983), 773-791. MR0711867
\bibitem[AF]{AF} Alves, C.O. and de Figueiredo, D.G.,
            \textit{Nonvariational elliptic systems},
             Discrete Contin. Dyn. Syst.
             \textbf{8} (2002), 289-302. MR1897684
\bibitem[B]{B} Bechah, A.,
            \textit{Positive solutions for a nonvariational
            quasilinear elliptic system in $\mathbb{R}^N$},
             Rev. R. Acad. Cienc. Exactas F¨ªs. Nat. (Esp.)
             \textbf{94} (2000), 1-7. MR1829496
\bibitem[BC]{BC} Br\'{e}zis, H. and Cabr\'{e}, X.,
            \textit{Some simple nonlinear PDE's without solutions},
             Boll. Unione Mat. Ital. Sez. B Artic. Ric. Mat. (8)
             \textbf{1} (1998), 223-262. MR1638143
\bibitem[BCMR]{BCMR} Br\'{e}zis, H., Cazenave, T., Martel, Y. and
             Ramiandrisoa,A.,
            \textit{Blow up for $u_t-\Delta u=g(u)$ revisited},
             Adv. Differential Equations
             \textbf{1} (1996), 73-90. MR1357955
\bibitem[BK]{BK} Br\'{e}zis, H. and Kato, T.,
            \textit{Remarks on the Schr\"{o}dinger operator with singular complex potentials},
             J. Math. Pures Appl. (9)
             \textbf{58} (1979), 137-151. MR0539217
\bibitem[BM]{BM} Birindelli, I. and Mitidieri, E.,
            \textit{Liouville theorems for elliptic inequalities and applications},
             Proc. Roy. Soc. Edinburgh
             \textbf{128A} (1998), 1217-1247. MR1664101
\bibitem[BuM]{BuM} Busca, J. and Manasevich, R.,
            \textit{A Liouville-type theorem for Lane-Emden system},
             Indiana Univ. Math. J.
             \textbf{51} (2002), 37-51. MR1896155
\bibitem[BT]{BT} Br\'{e}zis, H. and Turner, R.E.L.,
            \textit{On a class of superlinear elliptic problems},
             Comm. Partial Differential Equations
             \textbf{2} (1977), 601-614. MR0509489
\bibitem[BV]{BV} Bidaut-V\'{e}ron, M.-F. and Vivier, L.,
            \textit{An elliptic semilinear equation with source term involving boundary measures: the subcritical case},
             Rev. Mat. Iberoamericana
             \textbf{16} (2000), 477-513. MR1813326
\bibitem[C]{C} Cosner, C.,
            \textit{Positive solutions for superlinear elliptic systems without variational structure},
             Nonlinear Anal.
             \textbf{8} (1984), 1427-1436. MR0769404
\bibitem[Ch]{Ch} Chen, H.,
            \textit{Positive steady-state solutions of a non-linear reaction-diffusion system},
             Math. Methods Appl. Sci.
             \textbf{20} (1997), 625-634. MR1441724
\bibitem[CFM]{CFM} Clement, Ph., de Figueiredo, D.G. and Mitidieri, E.,
            \textit{Positive solutions of semilinear elliptic systems},
             Comm. Partial Differential Equations
             \textbf{17} (1992), 923-940. MR1177298
\bibitem[CFM$_2$]{CFM$_2$} Clement, Ph., de Figueiredo, D.G. and Mitidieri, E.,
            \textit{A priori estimates for positive solutions of semilinear elliptic systems via
                    Hardy-Sobolev inequalities},
             Pitman Res. Notes Math.
             \textbf{343} (1996), 73-91. MR1417272
\bibitem[CFMT]{CFMT} Clement, Ph., Fleckinger, J., Mitidieri, E. and de Th\'{e}lin,
               F.,
            \textit{Existence of positive solutions for a nonvariational quasilinear elliptic system},
             J. Differential Equations
             \textbf{166} (2000), 455-477. MR1781264
\bibitem[CMM]{CMM} Clement, Ph., Manasevich, R. and Mitidieri, E.,
            \textit{Positive solutions for a quasilinear system via blow up},
             Comm. Partial Differential Equations
             \textbf{18} (1993), 2071-2106. MR1249135
\bibitem[CFS]{CFS} Cuesta, M., de Figueiredo, D.G. and Srikant, P.N.,
            \textit{On a resonant-superlinear elliptic problem},
             Calc. Var. Partial Differential Equations
             \textbf{17} (2003), 221-233. MR1989831
\bibitem[DE]{DE} Dickstein, F. and Escobedo, M.,
             \textit{A maximum principle for semilinear parabolic systems and applications},
             Nonlinear Anal.
             \textbf{45} (2001), 825-837. MR1845028
\bibitem[DMP]{DMP}Del Pino, M., Musso, M. and Pacard, P.,
            \textit{Boundary singularities for weak solutions of semilinear elliptic problems},
             J. Funct. Anal.
             \textbf{253} (2007), 241--272. MR2362423
\bibitem[F]{F} de Figueiredo, D.G.,
            Semilinear elliptic systems. Nonlinear Functional Analysis and Applications to Differential Equations,
            Trieste 1997,World Sci. Publishing, River Edge, N.J.,
            1998, pp. 122-152.
\bibitem[FF]{FF} de Figueiredo, D.G. and Felmer, P.,
            \textit{A Liouville-type theorem for elliptic systems},
             Ann. Scuola Norm. Sup. Pisa Cl. Sci. (4)
             \textbf{21} (1994), 387-397. MR1310633
\bibitem[FLN]{FLN} de Figueiredo, D.G.,  Lions, P.-L. and Nussbaum, R.D. ,
            \textit{A priori estimates and existence of positive solutions of semilinear elliptic equations},
             J. Math. Pures Appl. (9)
             \textbf{61} (1982), 41-63. MR0664341
\bibitem[FSW]{FSW} Fila, M., Souplet, Ph. and Weissler, F.,
            \textit{Linear and nonlinear heat equations in $L^p_\delta$ spaces and universal bounds for
            global solutions},
             Math. Ann.
             \textbf{320} (2001), 87-113. MR1835063
\bibitem[FY]{FY} de Figueiredo, D.G. and Yang, J.,
            \textit{A priori bounds for positive solutions of a non-variational elliptic system},
             Comm. Partial Differential Equations
             \textbf{26} (2001), 2305-2321. MR1876419
\bibitem[GS]{GS} Gidas, B. and Spruck, J.,
            \textit{A priori bounds for positive solutions of nonlinear elliptic equations},
             Comm. Partial Differential Equations
             \textbf{6} (1981), 883-901. MR0619749
\bibitem[GW]{GW} Gu, Y.-G. and Wang, M.-X.,
            \textit{A semilinear parabolic system arising in the
                    nuclear reactors},
             Chinese Sci. Bull.
             \textbf{39} (1994), 1588-1592.
\bibitem[GW$_2$]{GW$_2$} Gu, Y.-G. and Wang, M.-X.,
            \textit{Existence of positive stationary solutions and threshold results for
                    a reaction-diffusion system},
             J. Differential Equations
             \textbf{130} (1996), 277-291. MR1410888
\bibitem[Li]{Li} Li, Y.-X., Liu, Q.-L. and Xie, Ch.-H.,
            \textit{Semilinear reaction diffusion systems of several component},
             J. Differential Equations
             \textbf{187} (2003), 510-519. MR1949453
\bibitem[Lou]{Lou} Lou, Y.,
            \textit{Necessary and sufficient condition for the existence
                   of positive solutions of certain cooperative system},
             Nonlinear Anal.
             \textbf{26} (1996), 1079-1095. MR1375651
\bibitem[JL]{JL}  Joseph, D.D. and Lundgren, T.S.,
            \textit{Quasilinear Dirichlet problems driven by positive sources},
             Arch. Ration. Mech. Anal.
             \textbf{49} (1972/73), 241-269. MR0340701
\bibitem[M]{M} Mitidieri, E.,
            \textit{Nonexistence of positive solutions of semilinear elliptic systems in $R\sp N$},
             Differential Integral Equations
             \textbf{9} (1996), 465-479. MR1371702
\bibitem[M$_2$]{M$_2$} Mitidieri, E.,
            \textit{A Rellich type identity and applications},
             Comm. Partial Differential Equations
             \textbf{18} (1993), 125-151. MR1211727
\bibitem[MR]{MR} McKenna, P.J. and Reichel, J.,
            \textit{A priori bounds for semilinear equations and
                    a new class of critical exponents for Lipschitz domains},
             J. Funct. Anal.
             \textbf{244} (2007), 220-246. MR2294482
\bibitem[NS]{NS} Ni, W.-M. and Sacks, P.,
            \textit{Singular behavior in nonlinear parabolic equations},
             Trans. Amer. Math. Soc.
             \textbf{287} (1985), 657-671. MR0768731
\bibitem[P]{P} Pacard, F.,
            \textit{Existence and convergence of positive weak solutions of
             $-\Delta u=u\sp {n/(n-2)}$ in bounded domains of $\bold R\sp n, n\geq 3$},
             Calc. Var. Partial Differential Equations
             \textbf{1} (1993), 243-265. MR1261546
\bibitem[PQS]{PQS}Polacik, P., Quittner, P. and Souplet, Ph.,
            \textit{Singularity and decay estimates in superlinear problems via Liouville-type
            theorems, I. Elliptic equations and systems},
             Duke Math. J.
             \textbf{139} (2007), 555-579. MR2350853
\bibitem[Q]{Q} Quittner, P.,
            \textit{Transition from decay to blow-up in a parabolic system},
             Arch. Math. (Brno)
             \textbf{34} (1998), 199-206. MR1629705
\bibitem[QS]{QS} Quittner, P. and Souplet, Ph.,
            \textit{A priori estimates and existence for elliptic systems via bootstrap in weighted Lebesgue spaces},
             Arch. Ration. Mech. Anal.
             \textbf{174} (2004), 49-81. MR2092996
\bibitem[QS$_2$]{QS$_2$} Quittner, P. and Souplet, Ph.,
            Superlinear parabolic problems: blow-up, global
            existence and steady states, Birkh\"{a}user Advanced
            Texts, Basel, Boston, Berlin, 2007.
\bibitem[RZ]{RZ} Reichel, W. and Zou, H.-H.,
            \textit{Non-existence results for semilinear cooperative elliptic systems via moving spheres},
             J. Differential Equations
             \textbf{161} (2000), 219-243. MR1740363
\bibitem[S]{S} Souplet, Ph.,
            \textit{Optimal regularity conditions for elliptic problems via
                    $L_\delta^p$-spaces},
             Duke Math. J.
             \textbf{127} (2005), 175-192. MR2126499
\bibitem[S$_2$]{S$_2$} Souplet, Ph.,
            \textit{A survey on $L^p_\delta$ spaces and their applications to nonlinear elliptic and parabolic problems},
             GAKUTO International Ser. Math. Sci. Appl.
             (Nonlinear Partial Differential Equations Their Appl.)
             \textbf{20} (2004), 464-479. MR2087491
\bibitem[St]{St} Struwe, M.,
            Variational methods. Applications to nonlinear partial
            differential equations and Hamiltonian systems,
            Springer, Berlin, 2000.
\bibitem[So]{So} Souto, M.,
            \textit{A priori estimate and existence of positive solutions
                    of nonlinear cooperative elliptic system},
             Differential Integral Equations
             \textbf{8} (1995), 1245-1258. MR1325555
\bibitem[SZ]{SZ} Serrin, J. and Zou, H.-H.,
            \textit{Non-existence of positive solutions of Lane-Emden systems},
             Differential Integral Equations
             \textbf{9} (1996), 635-653. MR1401429
\bibitem[SZ2]{SZ2} Serrin, J. and Zou, H.-H.,
            \textit{Existence of positive solutions of the Lane-Emden system},
             Atti Sem. Mat. Fis. Univ. Modena
             \textbf{46} suppl. (1998), 369-380. MR1645728
\bibitem[TV]{TV} de Th\'{e}lin, F. and V\'{e}lin, J.,
             \textit{Existence and nonexistence of nontrivial solutions for some nonlinear elliptic systems}
             Rev. Mat. Univ. Complut. Madrid
             \textbf{6} (1993), 153-194. MR1245030
\bibitem[Wang]{W} Wang, M.-X.,
             \textit{Global existence and finite time blow up for a reaction-diffusion system},
             Z. Angew. Math. Phys.
             \textbf{51} (2000), 160-167. MR1745297
\bibitem[Zou]{Zou} Zou, H.-H.,
            \textit{A priori estimates for a semilinear elliptic systems
                    without variational structure and their applications},
             Math. Ann.
             \textbf{323} (2002), 713-735. MR1924277
\bibitem[Zou$_2$]{Zou$_2$} Zou, H.-H.,
            \textit{A priori estimates and existence for strongly coupled semilinear cooperative elliptic systems},
             Comm. Partial Differential Equations
             \textbf{31} (2006), 735--773. MR2233039
\bibitem[Zh]{Zh} Zheng, S.-N.,
            \textit{Nonexistence of positive solutions to a
                    semilinear elliptic system and blowup estimates for a
                    reaction-diffusion system},
             J. Math. Anal. Appl.
             \textbf{232} (1999), 293-311. MR1683140
\bibitem[Zhao]{Zhao} Zhao, Z.-X.,
            \textit{Green function for Schr\"{o}dinger operator and conditioned Feynman-Kac gauge},
             J. Math. Anal. Appl.
             \textbf{116} (1986), 309-334. MR0842803
\bibitem[ZZ]{ZZ} Zhao, P.-H. and Zhong, C.-K.,
            \textit{On the infinitely many positive solutions of a supercritical elliptic problem},
            Nonlinear Anal.
            \textbf{44} (2001), 123-139. MR1815695
\end{thebibliography}
\end{document}